\documentclass{amsart}[11pt]
\usepackage{amsmath}
\usepackage{amssymb}
\usepackage{mathrsfs}
\usepackage[a4paper]{geometry}
\theoremstyle{plain}
\newtheorem{teo}{Theorem}[section]
\newtheorem{pro}{Proposition}[section]
\newtheorem{lem}{Lemma}[section]
\newtheorem{cor}{Corollary}[section]
\theoremstyle{remark}
\newtheorem*{nota}{Remark}
\newtheorem*{note}{Remarks}

\newcommand{\re}{\mathbb{R}}

\newcommand{\jac}{P_n^{r,s}}
\newcommand{\jam}{P_{\tau}^{r,s,\beta}}
\title{$\beta$-Jacobi processes}

\begin{document}
\maketitle
\centerline{N. Demni\footnote{SFB 701, Fakult\"at f\"ur Mathematik, universit\"at Bielefeld, Bielefeld, Germany. \\ email: demni@math.uni-bielefeld.de} }

\begin{abstract}
We define and study a $[0,1]^m$-valued process depending on three positive real parameters $p,q,\beta$ that specializes for $\beta=1,2$ to the eigenvalues process of the real and complex matrix Jacobi processes on the one hand and that  has the distribution of the $\beta$-Jacobi ensemble as stationary distribution on the other hand. 
We first prove that this process, called $\beta$-Jacobi process, is the unique strong solution of the stochastic differential equation defining it provided that $\beta > 0, p \wedge q > m - 1 + 1/\beta$. When specialized to $\beta=1,2$, our results actually improve well known results on eigenvalues of matrix Jacobi processes. While proving the strong uniqueness, the generator of the $\beta$-Jacobi process is mapped into  the radial part of the Dunkl-Cherednik Laplacian  associated with the non reduced root system of type $BC$. The transformed process is then valued in the principal Weyl alcove and this allows to define the Brownian motion in the Weyl alcove corresponding to all multiplicities equal one. 
Second, we determine, using stochastic calculus and a comparison theorem, the range of $\beta,p,q$ for which the $m$ components of the $\beta$-Jacobi process first collide, the smallest one reaches $0$ and the largest one reaches $1$. This is equivalent to the first hitting time of the boundary of the principal Weyl alcove by the transformed process. Finally, we write down its semi group density. 
\end{abstract}

\section{Introduction}

\subsection{Motivation}
The study of processes associated with root systems has known a considerable growth mainly after the birth of Dunkl operators which allowed to define Dunkl processes (\cite{Chy}). 
The projection of the latter processes on Weyl chambers associated with reduced root systems, known as radial Dunkl processes, revealed close connections to eigenvalues of some matrix-valued processes (\cite{Dem}) and to Brownian motions in Weyl chambers (\cite{De},\cite{Gra}, \cite{Karl}, \cite{Kat}). Indeed, the eigenvalues of symmetric and Hermitian Brownian matrices are radial Dunkl processes associated with type $A$ root systems while the square root of the eigenvalues of Wishart and Laguerre processes are radial Dunkl processes associated with type $B$ root systems. Besides, in the complex Hermitian setting, the eigenvalues process  is the $V$-process in Doob's sense of a multidimensional Brownian motion killed when it first hits the boundary of the Weyl chamber ($V$ stands for the Vandermonde function). 
Following the standard scheme in the theory of stochastic processes, one wonders how is the eigenvalues process of real and complex Jacobi matrix processes related to root systems and this was mainly behind our motivation. This connection is not new in its own and is implicitely written in \cite{Bee} where authors relate the multivariate Gauss hypergeometric function defined via Jack polynomials to the hypergeometric function associated with the root system of type $BC$. Loosely speaking, the eigenoperator of the firstly-mentioned function is the generator of a diffusion that we shall define and call the $\beta$-Jacobi process, where $\beta > 0$ is the inverse of the Jack parameter (see \cite{Bee} for more details and definitions). The name is justified by the fact that the $1$-and $2$ Jacobi processes fit the eigenvalues process of the real and complex Jacobi processes (\cite{Dou}) and also from the fact that the stationary distribution of the $\beta$-Jacobi process is given by the multivariate Beta distribution corresponding to the $\beta$-Jacobi ensemble (\cite{Che}).

\subsection{From matrix Jacobi processes to $\beta$-Jacobi processes}
The real (complex) matrix Jacobi process was introduced and studied in \cite{Dou}. It was defined as the radial part of a left corner of a $n \times n$ orthogonal (unitary) Brownian motion. The latter, say $\Theta_n$, may be defined via the heat kernel in the orthogonal (unitary) group (\cite{Bia}), and by a left corner of size $m \times p$ we mean the compression of $\Theta_n$ by two projections $P_m,Q_p$ of ranks $m \leq n,p \leq n $, that is $P_m\Theta_n Q_p$. The real (complex) matrix Jacobi process of parameters $(p,q = n-p)$ is then defined by $J_m:= P_m\Theta_nQ_p\Theta_n^{\star}P_m$ where $\Theta_n^{\star}$ stands for the adjoint of $\Theta_n$. This definition is extended to real parameters $p,q$ via stochastic differential equations. Let $\beta = 1,2$ respectively , then when $p\wedge q \geq (m-1)+ 2/\beta$, the eigenvalues process is the unique strong solution of (\cite{Dou} p. 135):
{\small \begin{equation}\label{JE}
d\lambda_i(t) = 2\sqrt{(\lambda_i(t)(1-\lambda_i(t))}d\nu_i(t) + \beta \left[(p - (p+q)\lambda_i(t)) + 
 \sum_{j \neq i}\frac{\lambda_i(t)(1-\lambda_j(t)) + \lambda_j(t)(1-\lambda_i(t))}{\lambda_i(t) - \lambda_j(t)}\right] dt
\end{equation}}
where $(\nu_i)_{i=1}^m$ are independent Brownian motions and $0 < \lambda_m \leq \cdots \leq  \lambda_1 < 1$. In fact, the eigenvalues remain in $]0,1[$ when 
$p \wedge q \geq (m-1) + 2/\beta$. When one considers a Haar orthogonal (unitary) matrix instead of an orthogonal (unitary) Brownian motion, the distribution of $J$ no more depends on time (stationary) and its eigenvalues are distributed at any time $t$ according to the multivariate Beta distribution: 
\begin{equation*}
Z_m^{p,q,\beta}(\lambda) := C_{n,p,q, \beta} \prod_{i=1}^m\lambda_i^{p-\beta(m-1)/2}(1-\lambda_i)^{q - \beta(m-1)/2} \prod_{1\leq i < j\leq m}|\lambda_i - \lambda_j|^{\beta}
\end{equation*}   
for $\beta=1,2$ respectively and for some normalizing constant $C_{n,p,q, \beta}$. In the same way the $\beta$-Jacobi ensemble (\cite{Kil}) generalizes the JOE, JUE, JSE (\cite{Che}), we introduce the $\beta$-Jacobi process as a solution, whenever it exists, of (\ref{JE}) for arbitrary $\beta > 0$. We will prove that (\ref{JE}) has a unique strong solution for $p \wedge q > (m-1)+1/\beta$, for any time $t$ and for any starting point $0 \leq \lambda_m(0) \leq \cdots \leq \lambda_1(0) \leq 1$. This result actually improves the strong existence and uniqueness result derived in \cite{Dou}. To proceed, we will instead focus on the process $\phi$ defined by $\phi_i := \arcsin\sqrt{\lambda_i},\,1 \leq i \leq m$ so that $0 \leq \phi_m \leq \dots \leq \phi_1 \leq \pi/2$. The process $\phi$ is shown to satisfy 
\begin{equation*}
d\phi(t) = d\nu(t) + \nabla \sum_{\alpha \in R_+}k(\alpha)\log(\sin\langle \alpha,\phi(t) \rangle)dt = d\nu(t) + \sum_{\alpha \in R_+}k(\alpha)\cot(\langle \alpha,\phi(t) \rangle) dt \,\alpha,
\end{equation*}
where $\nu = (\nu_i)_{i=1}^m$, $\langle \cdot \rangle$ is the usual inner product in $\re^m$, $R_+$ is a positive system of the non reduced root system of type $BC$ and 
$k: \alpha \mapsto k(\alpha)$ is a positive multiplicity function, that is, $k(\alpha) = k(\gamma)$ if and only if $\alpha,\gamma$ are conjugated under the action of the reflections group (\cite{Hum}). The strong uniqueness will follow from a scheme developed in \cite{Dem} in order to prove a similar result for radial Dunkl processes. However, this setting needs more care since the state space of $\phi$ (or $\lambda$) is a bounded domain while Weyl chambers are not. In fact, $\phi/\pi$ is valued in the principal Weyl alcove $A_{\circ}$ of type $BC$ 
(\cite{Hum}, p. 89) therefore one defines the Brownian motion in the Weyl alcove of type $BC$ as the process $\phi/\pi$ corresponding to $k(\alpha) \equiv 1$. Unlike Brownian motions in Weyl chambers and the Brownian motion in the Weyl alcove of type $A$ (see below), this process is not a $h$-process of a $m$-dimensional Brownian motion killed when it first reaches the boundary $\partial A_{\circ}$ of $A_{\circ}$, where
\begin{equation*}
h(\phi) := \prod_{\alpha \in R_+}\sin(<\alpha,\phi>).
\end{equation*}
More precisely, $h$ is shown to satisfy $(\Delta/2) h > ch, c < 0$. The end of the paper is devoted to the semi group density of $\lambda$ expressed as a bilinear series of multivariate Jacobi polynomials (\cite{Lass1}).  

\section{Strong uniqueness}    
\subsection{The transformed process} 
As mentioned before, the connection to root systems was implicitely written in \cite{Bee} and was done via differential operators which are the generators of $\lambda$ and $\phi$. 
For the reader's convenience, we display the details and derive the stochastic differential equation satisfied by $\phi$. 
Let $0 < \lambda_m < \cdots < \lambda_1 < 1$ and consider (\ref{JE}) up to 
\begin{equation*}
\inf\{t, \, \lambda_m(t) = 0\} \wedge \inf\{t,\, \lambda_1(t) = 1\} \wedge\{t, \, \lambda_i(t) = \lambda_j(t) \, \textrm{for some}\, (i,j)\}.
\end{equation*}
Let $\phi_i  := \arcsin(\sqrt\lambda_i) := s(\lambda_i), 1\leq i \leq m$, then $0< \phi_m(0) < \dots <\phi_1(0) < \pi/2$ and:
\begin{equation*}
s'(\lambda_i) = \frac{1}{\sin 2\phi_i},\qquad s''(\lambda_i) = \frac{2(2\sin^2\phi_i - 1)}{\sin^3 2\phi_i} = -\frac{2\cos 2\phi_i}{\sin^3 2\phi_i}.
\end{equation*} 
Using
\begin{equation*}
\sin^2\phi_i - \sin^2 \phi_j = \sin(\phi_i + \phi_j)  \sin(\phi_i - \phi_j),
\end{equation*}
\begin{equation*}
\sin^2\phi_i\cos^2\phi_j + \cos^2\phi_i\sin^2\phi_j = \frac{1}{2}[\sin^2(\phi_i+\phi_j) + \sin^2(\phi_i-\phi_j)],
\end{equation*}
then It\^o's formula gives: 
\begin{align*}
d\phi_i(t) & = d\nu_i(t) + \left[\beta\frac{(p - (p+q)\sin^2\phi_i)}{\sin 2\phi_i} - \cot 2\phi_i\right]dt 
 + \frac{\beta}{2}\frac{dt}{\sin2\phi_i(t)}\sum_{j \neq i}\frac{\sin^2(\phi_i(t)+\phi_j(t)) + \sin^2(\phi_i(t)-\phi_j(t))}{\sin(\phi_i(t) + \phi_j(t)) \sin(\phi_i(t) - \phi_j(t))}.
\end{align*}
Using 
\begin{equation*}
\sin^2\phi_i = \frac{1-\cos 2\phi_i}{2},\quad  \frac{1}{\sin 2\phi_i} = \cot \phi_i - \cot2\phi_i, 
\end{equation*}
then
\begin{align*}
\beta\frac{(p - (p+q)\sin^2\phi_i)}{\sin 2\phi_i} - \cot 2\phi_i  = \beta\frac{(p-q)}{2}\cot(\phi_i(t)) + [\beta q -1]\cot (2\phi_i(t)). 
\end{align*}
Using
\begin{eqnarray*}
\sin2\phi_i &=& [\cot(\phi_i+\phi_j) + \cot(\phi_i-\phi_j)]\sin(\phi_i+\phi_j)\sin(\phi_i-\phi_j),\\
\frac{1}{\sin^2z} &=& 1+\cot^2 z,
\end{eqnarray*} 
one gets 
\begin{align*}
&\frac{1}{2\sin2\phi_i(t)}\frac{\sin^2(\phi_i(t)+\phi_j(t)) + \sin^2(\phi_i(t)-\phi_j(t))}{\sin(\phi_i(t) + \phi_j(t)) \sin(\phi_i(t) - \phi_j(t))} = 
\frac{1}{2}\frac{[1/\sin^2(\phi_i(t)+\phi_j(t))] + [1/\sin^2(\phi_i(t)-\phi_j(t))]}{\cot(\phi_i(t)+\phi_j(t)) + \cot(\phi_i(t)-\phi_j(t))}\\&
= \frac{1}{2}\frac{\cot^2(\phi_i(t)+\phi_j(t)) + \cot^2(\phi_i(t)+\phi_j(t)) + 2}{\cot(\phi_i(t)+\phi_j(t)) + \cot(\phi_i(t)-\phi_j(t))}
\\& = \frac{1-\cot(\phi_i(t)+\phi_j(t))\cot(\phi_i(t)-\phi_j(t))}{\cot(\phi_i(t)+\phi_j(t)) + \cot(\phi_i(t)-\phi_j(t))}+\frac{\cot(\phi_i(t)+\phi_j(t)) + \cot(\phi_i(t) - \phi_j(t))}{2}.
\end{align*}
Finally, using
\begin{equation*}
\cot(u+v) = \frac{\cot(u) \cot(v)-1}{\cot(u) + \cot(v)},
\end{equation*}
one gets 
\begin{align}\label{JE1}
d\phi_i(t) = d\nu_i(t)  +\left[k_0\cot \phi_i + k_1 \cot 2\phi_i(t)dt +  k_2\sum_{i\neq j}[\cot(\phi_i+\phi_j) + \cot(\phi_i -\phi_j)]\right]dt
\end{align}where 
\begin{equation}\label{MF}
2k_0 = \beta(p-q),\quad k_1 = \beta(q-(m-1)) -1, \quad 2k_2 = \beta.
\end{equation}
Easy computations show that $\pi/2 - \phi$ satisfies (\ref{JE1}) with $(p,q)$ intertwined.

\subsection{Existence and uniqueness of a strong solution.}
From (\ref{JE1}), the generator of $\phi$ acts on smooth functions as (it is the trigonometric version of the $W$-invariant Dunkl-Cherednik Laplacian, \cite{Scha}): 
\begin{equation*}
\mathscr{L} := \frac{1}{2}\Delta - \langle\nabla, \nabla\Phi(\cdot)\rangle, \quad \Phi(\phi) := -\sum_{\alpha \in R_+}k(\alpha)\log\sin(\langle\alpha,\phi\rangle),
\end{equation*}
where   
\begin{eqnarray*}
R &=& \{\pm e_i, \, \pm 2e_i,\, 1 \leq i \leq m,\, \pm(e_i \pm e_j),\, 1 \leq i < j \leq m\}, \\
R_+ & = &  \{e_i, \,  2e_i,\, 1 \leq i \leq m,\, (e_i \pm e_j),\, 1 \leq i < j \leq m\}, \\
S &=& \{e_i - e_{i+1}, \,1 \leq i \leq m-1,\, e_m\},
\end{eqnarray*}  
 $(e_i)_{i=1}^m$ is the canonical basis of $\re^m$ and $S$ is the simple system corresponding to the positive system $R_+$ (\cite{Hum}). $R$ is known as the non reduced root system of type $BC$ (\cite{Hum} p. 41) and the action of the reflections group on $R$ gives rise to three orbits so that the multiplicity function takes three values corresponding to 
 $\{\pm e_i,\, 1 \leq i \leq m\}, \{\pm 2e_i,\, 1 \leq i \leq m\}, \{\pm(e_i \pm e_j),\, 1 \leq i < j \leq m\}$. With regard to (\ref{JE1}), 
 $k(e_i) = k_0, 2k(2e_i) = k_1, k(e_i \pm e_j) = k_2$. Setting $\tilde{\phi}_i := \phi_i/\pi$, then $\tilde{\phi}$ is valued in the closure of the {\it principal Weyl alcove} (\cite{Hum}, p. 90) defined by : 
\begin{equation*}
A_{\circ} = \{\tilde{\phi} \in \re^m, \, \langle \alpha,\tilde{\phi}\rangle > 0 \, \textrm{for all}\, \alpha \in S, \, \langle \tilde{\alpha},\tilde{\phi}\rangle  < 1\} 
= \{\tilde{\phi} \in \re^m, \,0 < \tilde{\phi}_m < \cdots < \tilde{\phi_1} < 1/2\}
\end{equation*}
where $\tilde{\alpha} = 2e_1$ is the unique highest root (\cite{Hum} p. 40). When $k_0 = 0  \Leftrightarrow p=q$, (\ref{JE1}) invovles the root system of type $C$ given by the data (\cite{Hum} p. 42)  
\begin{eqnarray*}
R&=& \{\pm e_i \pm e_j,\,1 \leq i < j \leq m,\, \pm 2e_i, \,1 \leq i \leq m\}\\
R_+&=& \{e_i \pm e_j,\,1 \leq i < j \leq m,\, 2e_i, \,1 \leq i \leq m\}\\
S &=& \{e_i - e_{i+1}, \,1 \leq i \leq m-1,\, e_m\}
\end{eqnarray*} 
and we refer to the $\beta$-Jacobi process as the $\beta$-ultraspherical process. Now, we proceed to the proof of the strong uniqueness for (\ref{JE1}) subject to $k_0, k_1, k_2 > 0$. Once we did, we apply this result to $\pi/2-\phi$ and to the ultraspherical $\beta$-Jacobi process to deduce the strong uniqueness for $p \wedge q > (m-1)+1/\beta$. 
For, we need to briefly recall the scheme developed in (\cite{Dem}) used to deal with radial Dunkl processes and inspired from \cite{Cepa}, \cite{Cepa1}. 
First, a result from \cite{Cepa}, \cite{Cepa1} states that  
\begin{equation*}
d\phi(t) = d\nu(t) - \nabla \Phi(\phi(t)) + n(\phi(t))dL_t 
\end{equation*} 
where $n(\phi(t))$ is a unitary inward normal vector to $\overline{A_{\circ}}$ at $\phi(t) \in \partial A_{\circ}$, that is (\cite{Cepa1})  
\begin{equation}\label{INV}
\langle \phi(t) - a, n(\phi(t)) \rangle \, \leq 0, \quad \forall a \in \overline{A_{\circ}}, 
\end{equation}  
and $L$ is the boundary process satisfying: 
\begin{equation*} 
dL_t = {\bf 1}_{\{\phi(t) \in \partial A_{\circ}\}} dL_t, 
\end{equation*} 
has a unique strong solution for all $t \geq 0$ and $\phi(0) \in \overline{A_{\circ}}$. Second, it remains to prove that $L$ vanishes. For, we need to prove two Lemmas: 
\begin{lem}\label{L1}
Let 
\begin{eqnarray*}
H_{\alpha} &:=& \{\phi \in \re^m,\, \langle \alpha,\phi \rangle = 0\},\, \alpha \in S, \\
H_{\tilde{\alpha}} &:=& \{\phi \in \re^m, \, \langle \tilde{\alpha},\phi\rangle = \pi\} 
\end{eqnarray*}
be the walls of $\overline{A_{\circ}}$. Then
\begin{equation*}
\partial A_{\circ}  = \cup_{\alpha \in S \cup \{\tilde{\alpha}\}} H_{\alpha} \cap \{\phi,\,\langle n(\phi),\alpha\rangle \neq 0\},  
\end{equation*}
that is if $\phi \in \partial A_{\circ}$, then there exists $\alpha \in S \cup \{\tilde{\alpha}\}$ such that  $\phi \in H_{\alpha}$ and $\langle n(\phi),\alpha\rangle \neq 0$.
\end{lem}
Once we prove this Lemma, the following trivially holds 
 \begin{equation*}
L_t  = {\bf 1}_{\{\phi(t) \in \partial A_{\circ}\}} dL_t \leq \sum_{\alpha \in S \cup \{\tilde{\alpha}\}} {\bf 1}_{\{\phi(t) \in H_{\alpha}\}} {\bf 1}_{\{\langle n(\phi(t)),\alpha\rangle \neq 0\}}dL_t.
\end{equation*}
As a matter of fact, the strong existence and uniqueness result will follow after proving that
\begin{lem}\label{L2}
\begin{equation*} 
{\bf 1}_{\{\phi(t) \in H_{\alpha}\}} \langle n(\phi(t), dL_t\rangle = 0, \quad \alpha \in S \cup \{\tilde{\alpha}\}. 
\end{equation*} 
\end{lem}
{\it Proof of Lemma \ref{L1}}: let $\phi \in \partial A_{\circ}$ then either $\phi$  belongs to $H_{\tilde{\alpha}}$ and one claims that $n(\phi) = -e_1 = -\tilde{\alpha}/2$ is a unitary inward normal vector at $\phi$, therefore $\langle n(\phi),\tilde{\alpha}\rangle \neq 0$. This claim follows immediately from (\ref{INV}) and from the fact that $\phi_1 = \pi/2$. 
Or ($\phi \in \partial A_{\circ} \setminus H_{\tilde{\alpha}}$) there only exists $\alpha \in S$ such that $\phi \in H_{\alpha}$ and one proceeds as in Lemma 2 in \cite{Dem} but with more care since the Weyl alcove is a bounded domain while the Weyl chamber is not. Hence, assume that for all $\alpha \in S$ such that $\langle \alpha,\phi \rangle = 0$, one has $\langle n(\phi), \alpha\rangle = 0$.  The idea is to find $\epsilon > 0$ such that $\phi- \epsilon n(\phi) \in \overline{A_{\circ}}$ and then conclude that $n(\phi) = 0$ by substituting $a =  \phi- \epsilon n(\phi)$ in (\ref{INV}). \\
If $\langle n(\phi), \alpha\rangle = 0$ for all $\alpha \in S$, then $\langle n(\phi),\tilde{\alpha}\rangle = 0$ and $\phi- \epsilon n(\phi) \in \overline{A_{\circ}}$ for all $\epsilon > 0$ and 
$\phi \in \partial A_{\circ}$. If $\langle n(\phi), \alpha\rangle \neq 0$ for some $\alpha \in S$, then $\langle \alpha,\phi \rangle > 0$. For those simple roots, either $\langle \alpha,n(\phi) \rangle < 0$ thus $\langle \tilde{\alpha}, n(\phi) \rangle < 0$ ($\tilde{\alpha} \in R_+$) and one seeks $\epsilon$ such that 
\begin{equation*}
\langle \tilde{\alpha},\phi - \epsilon n(\phi) \rangle \leq 1.
\end{equation*}
Such an $\epsilon$ exists since we assumed that $\langle \tilde{\alpha},\phi \rangle < 1$. Finally, if for those simple roots, some are such that $\langle \alpha,n(\phi) \rangle > 0$,
then one seeks $\epsilon$ such that 
\begin{equation*}
0 < \epsilon < \frac{\langle \phi,\alpha \rangle}{\langle n(\phi),\alpha \rangle}
\end{equation*}
for all $\alpha$ such that $\langle \phi,\alpha \rangle > 0, \langle n(\phi),\alpha \rangle > 0$ and such that 
$-\epsilon \langle \tilde{\alpha}, n(\phi) \rangle \leq 1 - \langle \tilde{\alpha},\phi\rangle$. The first condition is satisfied provided that 
\begin{equation*}
0< \epsilon < \min_{\langle \phi,\alpha \rangle > 0, \langle n(\phi),\alpha \rangle > 0} \frac{\langle \phi,\alpha \rangle}{\langle n(\phi),\alpha \rangle},
\end{equation*}
while for the second we need to write $\tilde{\alpha} = \sum_{\alpha \in S}a_{\alpha}\alpha, a_{\alpha} \geq 0$ so that 
\begin{equation*}
\langle \tilde{\alpha},n(\phi)\rangle = \sum_{\langle \alpha,n(\phi)\rangle > 0} a_{\alpha} \langle\alpha,n(\phi)\rangle +  \sum_{\langle \alpha,n(\phi)\rangle < 0} a_{\alpha} \langle\alpha,n(\phi)\rangle \geq \sum_{\langle \alpha,n(\phi)\rangle < 0} a_{\alpha} \langle\alpha,n(\phi)\rangle.
\end{equation*} 
Thus, 
\begin{equation*}
-\epsilon \langle \tilde{\alpha},n(\phi)\rangle \leq -\epsilon \sum_{\langle \alpha,n(\phi)\rangle < 0} a_{\alpha} \langle\alpha,n(\phi)\rangle .
\end{equation*}
If there is no $\alpha$ such that $\langle \phi,\alpha \rangle > 0, \langle n(\phi),\alpha \rangle < 0$, then $\langle \tilde{\alpha}, n(\phi) \rangle > 0$ since $\tilde{\alpha} \in R_+$ and
\begin{equation*}
\langle \tilde{\alpha},\phi - \epsilon n(\phi)  \rangle < \langle \tilde{\alpha}, \phi \rangle < 1.
\end{equation*}
Otherwise 
\begin{equation*}
0 < \epsilon <   [- \sum_{\langle \alpha,n(\phi)\rangle < 0} a_{\alpha} \langle\alpha,n(\phi)\rangle]^{-1}(1 - \langle \tilde{\alpha},\phi\rangle)
\end{equation*}
and one chooses the smallest $\epsilon$ which satisfies both conditions. The Lemma is proved. $\hfill \blacksquare$\\

{\it Proof of Lemma \ref{L2}}: when $\alpha \in S$, the proof is exactly the same as in Lemma 1 in \cite{Dem} when substituting $x \mapsto -\ln \langle \alpha,x\rangle$ by 
$\phi \mapsto -\log \sin \langle \alpha,\phi \rangle$. Thus, one has to deal with the additional term corresponding to $H_{\tilde{\alpha}}$.
Let $\theta(\phi) := -\log \sin (\phi)$,  then the occupation density formula yields (\cite{Rev}):
\begin{align*}
\int_0^{\pi}L_t^a(\pi - \langle\tilde{\alpha},\phi\rangle)|\theta^{'}(a)|da &= \langle \tilde{\alpha},\tilde{\alpha}\rangle  \int_0^t |\theta^{'}(\pi - \langle\tilde{\alpha},\phi(s)\rangle)| ds 
= 4  \int_0^t |\theta^{'}(\langle \tilde{\alpha},\phi(s)\rangle)| ds 
\end{align*}
since $\cot(\pi-z) = -\cot(z)$, where $L_t^a(\pi - \langle\tilde{\alpha},\phi\rangle)$ is the local time at time $t$ and at level $a$ of the real-valued semimartingale 
$\pi - \langle \tilde{\alpha},\phi\rangle$. Then, following line by line the proof of Lemma 1 in \cite{Dem}, we get
\begin{equation*}
\int_0^{\pi}L_t^a(\pi - \langle\tilde{\alpha},\phi\rangle)|\theta^{'}(a)|da < \infty
\end{equation*}
which implies that $L_t^0(\pi - \langle\tilde{\alpha},\phi\rangle) = 0$ for all $t$ since $\theta'$ blows up at $0$. Then, we use Tanaka's formula to compute 
\begin{equation*}
0 = d[\pi - \langle \tilde{\alpha},\phi(t)\rangle - (\pi - \langle \tilde{\alpha},\phi(t)\rangle)^+] = - {\bf 1}_{\{\pi - \tilde{\alpha},\phi(t)\rangle = 0\}} d \langle \tilde{\alpha},\phi(t)\rangle 
\end{equation*}
and finally use similar arguments to prove the statement of the Lemma. $\hfill \blacksquare$ \\
Thus we proved that  
\begin{teo} 
If $k_0 > 0,\,k_1 > 0,\, k_2 > 0$, then (\ref{JE1}) has a unique strong solution for all $t > 0$ and $\phi(0) \in \overline{A_{\circ}}$.
\end{teo}
Since $\pi/2 - \phi$ satisfies (\ref{JE1}) with $p,q$ intertwined and since the case $p=q \Leftrightarrow k_0= 0$ involves the reduced root system of type $C$ which share the same Weyl alcove and the same highest root with the root system of type $BC$, our proof applies to those cases. We finally deduce that 
\begin{cor} 
If $p \wedge q > (m-1) + 1/\beta, \beta > 0$, then (\ref{JE1}) has a unique strong solution for all $t > 0$ and $\phi(0) \in \overline{A_{\circ}}$.
\end{cor}
\begin{note}
1/Our result simplifies to $p \wedge q > m$ in the real case $\beta = 1$ and $p \wedge q > m-1/2$ in the complex one $\beta =2$. \\
2/The process $\phi/\pi$ may be interpreted as particles in the interval $[0,1/2]$. For particles on the circle, see \cite{Cepa2}, \cite{Hob}. In that case, $\phi_m < \cdots < \phi_1 < \phi_m +2\pi$ so that $\phi/2\pi$ is valued in the Weyl alcove of type $A$ (the highest root is $\phi_1-\phi_m$, \cite{Hum}) p. 41).
\end{note}

\subsection{Brownian motion in Weyl alcoves}
It is known that Brownian motions in Weyl chambers are radial Dunkl processes with a multiplicity function that equals one and that they are the unique strong solution of some stochastic differential equation with a rational singular drift for all $t$ (\cite{Dem},\cite{Gra}). Similarly, the Brownian motion in the principal Weyl alcove of type $A$ or equivalently the eigenvalues process of a unitary Brownian motion satisfies 
(see the end of \cite{Hob}) 
\begin{equation*}
d\phi_i(t) = d\nu_i(t) + \sum_{j \neq i}\cot(\phi_i(t)-\phi_j(t))dt,
\end{equation*} 
which also has a unique strong solution for all time $t$ (\cite{Cepa2}). As a matter of fact, we define the Brownian motion in the principal Weyl alcove of type $BC$ as the unique strong solution of (\ref{JE1}) with $k(\alpha) = 1$ for all $\alpha \in R$, that is $k_0 = k_2 = 1, k_1 = 2$ or $\beta = 2, q = m+1/2, p = q+1 = m+3/2$. Its generator reads 
\begin{equation*}
\mathscr{L} = \frac{1}{2}\Delta + \sum_{\alpha \in R_+}\langle \nabla, \nabla \log\sin(\langle\alpha,\phi\rangle) := \frac{1}{2}\Delta + \langle \nabla, \nabla \log h(\phi)\rangle,
\end{equation*}
where 
\begin{equation*}
h(\phi) := \prod_{\alpha \in R_+}\sin(\langle\alpha,\phi \rangle) = \prod_{j=1}^m\sin(\phi_j)\sin(2\phi_j)\prod_{j< l}\sin(\phi_j- \phi_l)\sin(\phi_j+\phi_l). 
\end{equation*}
Then, $h$ is strictly positive on $A_{\circ}$ and vanishes for $\phi \in \partial A_{\circ}$. Unfortunately, unlike the case of Brownian motions in Weyl chambers or the case of the Brownian motion in the principal Weyl alcove of type $A$, this process is not a $h$-process of a mutidimensional Brownian motion. Indeed, we can prove (see Appendix) that $h$ is not an eigenfunction of $\Delta/2$ and rather satisfies $\Delta h > ch$ for some strictly negative constant $c$. Nevertheless, when $\beta=2$, $\lambda$ is known to be the $V$-transform ($V$ stands for the Vandermonde function) in Doob's sense of $m$ independent real Jacobi processes of special parameters  conditioned never to collide (\cite{Dou}, p. 141). 

\section{The first hitting time of $\partial A_{\circ}$}
We define  the first hitting time of $\partial A_{\circ}$ by 
\begin{equation*}
T_0 = \inf\{t > 0,\, (\phi(t)/\pi) \in \partial A_{\circ}\} = T_{\tilde{\alpha}} \wedge \inf_{\alpha \in S} T_{\alpha}
\end{equation*} 
where 
\begin{eqnarray*}
T_{\alpha} &:=& \inf\{t > 0,\, <\alpha,\phi(t)> = 0\},\\
T_{\tilde{\alpha}} &:= & \inf\{t > 0, <\tilde{\alpha}, \phi(t)> = 2\phi_1 = \pi\}, 
\end{eqnarray*}
and $\phi$ is the unique strong solution for all $t \geq 0$ of (\ref{JE1}) with $R=BC_m$ and $p \wedge q > (m-1) + 1/\beta$. We will prove via stochastic calculus that 
\begin{pro}
\begin{itemize}\noindent
\item[$\bullet$] If $0 < k_2 < 1/2 \Leftrightarrow 0 < \beta < 1$, then $T_{\alpha} < \infty$ a.s. for $\alpha \in \{e_i - e_{i+1}, 1 \leq i \leq m-1\}$.
\item[$\bullet$] If $0< k_0 + k_1 +1< 2 \Leftrightarrow 0 < p - (m-1) < 2/\beta $, then $T_{e_m} < \infty$ a.s.. 
\item[$\bullet$] $0 < q - (m-1) < 2/\beta$, then  $T_{\tilde{\alpha}} = T_{2e_1} < \infty$ a.s..
\end{itemize}
\end{pro}
Recall that for a radial Dunkl process $X$ associated with a reduced root system $R$ and a Weyl chamber $C$, we compared the one dimensional process $\langle \alpha,X\rangle, \alpha \in S$ to a Bessel process in order to prove the finiteness of the first hitting time of $\partial C$ (\cite{Dem}). Here, we shall compare $\langle \alpha,\phi \rangle, \alpha \in S$ to a real Jacobi process and use that $\pi/2 - \phi$ satisfies (\ref{JE1}) to deal with $T_{\tilde{\alpha}}$. However, the case we have in hands is more delicate since the root system of type $BC$ is non reduced. We start with recalling some needed facts about Jacobi processes from \cite{Dou}, \cite{War}. 
  
\subsection{On real Jacobi processes} 
The real Jacobi process of parameters $d,d' \geq 0$ is the unique strong solution  of 
\begin{equation*}
dJ_t = 2\sqrt{J_t(1-J_t)}d\gamma_t + (d - (d+d')J_t) dt, \quad d,d' \geq 0, 
\end{equation*}
where $\gamma$ is a real Brownian motion. It is known from (\cite{Dou}) that $J$ remains in $]0,1[$ when $d \wedge d' \geq 2$, hits $0$ if $0 \leq d < 2$ and hits $1$ if $0 \leq d' < 2$. 
According to our result on strong uniqueness specialized to the rank one case $BC_1$, $\psi := \arcsin \sqrt{J}$ is the unique strong solution of 
\begin{equation}\label{JE2}
d\psi_t = d\gamma_t + [(d-d')\cot(\psi_t) + (d'-1)\cot(2\psi_t)]dt.
\end{equation}
subject to $d \wedge d' > 1$.

\subsection{Proof of the Proposition} Let $\alpha_0 \in R_+$, then from (\ref{JE1}), one has
\begin{align*}
d<\alpha_0,\phi(t)> = ||\alpha_0||d\gamma_t + k_2||\alpha_0||^2\cot\langle \alpha_0,\phi(t)\rangle dt+ 
\sum_{\alpha \in R_+\setminus \alpha_0}k(\alpha)a(\alpha)\cot \langle \alpha,\phi(t)\rangle dt
\end{align*} 
where $a(\alpha) = <\alpha_0,\alpha>$. Now, let $\alpha_0 \in \{e_i - e_{i+1}, \, 1 \leq i \leq m-1\}$ and let $p \geq q$. Denote $\sigma_0$ the reflection with respect to $H_{\alpha_0}$. Then, for all $\alpha \in R_+ \setminus \alpha_0$, easy computations show that $\langle \sigma_0(\alpha),\alpha_0\rangle = -a(\alpha)$. Since $\sigma_0(\alpha)$ belongs to the same orbit containing $\alpha$ and is positive (see \cite{Hum} p. 10), one has 
\begin{equation*}
F_t:= \sum_{\alpha \in R_+\setminus  \alpha_0}k(\alpha)a(\alpha)\cot\langle \alpha,\phi(t)\rangle = 
\sum_{\substack{\alpha \in R_+\setminus \alpha_0 \\ a(\alpha) > 0}} k(\alpha)a(\alpha)[\cot\langle \alpha,\phi(t)\rangle - \cot\langle \sigma_0(\alpha),\phi(t)\rangle].
\end{equation*}    
 
As a result, one writes for any $t \geq 0$: 
\begin{align*}
d(\phi_i(t) - \phi_{i+1}(t)) = \sqrt{2}d\gamma_t + 2k_2 \cot(\phi_i(t) - \phi_{i+1}(t)) dt + F_t dt.
\end{align*}  
This drift is strictly negative on $\{T_{\alpha_0} = \infty\}$ since $\phi \mapsto \cot \phi$ is a decreasing function and since $\langle \alpha_0, \phi(t)\rangle > 0$ so that: 
\begin{equation*}
\langle \alpha - \sigma_0(\alpha),\phi(t)\rangle = 2\frac{a(\alpha)}{||\alpha||^2}\langle \alpha_0,\phi(t)\rangle > 0 .
\end{equation*}
Using Proposition 2. 18. p. 293 and Exercice 2. 19. p. 294 in \cite{Kar} , one gets 
\begin{equation*}
\mathbb{P}_{\phi(0)}(\forall t \geq 0,\, \langle \alpha_0,\phi(t)\rangle \leq Z_{t}) = 1
\end{equation*} 
where 
\begin{equation*}
dZ_t = \sqrt{2}d\gamma_t + 2k_2\cot(Z_t) dt, \quad Z_0 = \langle \alpha_0,\phi(0)\rangle 
\end{equation*}
is defined on the same probability space as $\phi$. According to (\ref{JE2}), $(\sin^2Z_{t/2})_{t \geq 0}$ is a real Jacobi process of parameters $d = 2k_2+1, d'=1$ and the first statement of the Proposition is proved when $p \geq q$. Applying the same scheme to $\pi/2-\phi$ when $q \geq p$, we get the desired first statement. \\ 
Now, let $\alpha_0 = e_m$. Compared to the previous case, one has to take $\{e_m,2e_m\}$ out of the drift term of $\langle e_m,\phi \rangle$ to get that $\sigma_0(\alpha) \in R_+$ for $\alpha \in R_+ \setminus \{e_m,2e_m\}$. The last fact is easily checked since for $\alpha = e_i \pm e_j$ this amounts to consider the reduced root system of type $B$ (\cite{Hum} p. 10), otherwise $e_m$ is orthogonal to $\{e_i, 2e_i,\,i \neq m\}$ so that $\sigma_0$ fixes pointwise this set. Accordingly: 
\begin{align*}
d<\alpha_0,\phi(t)> = d\phi_m(t) = d\gamma_t + k_0\cot(\phi_m(t)) dt + k_1\cot(2\phi_m(t))+F_t dt
\end{align*}  where 
\begin{align*}
F_t &= \sum_{\substack{\alpha \in R_+^1\setminus \{e_m,2e_m\} \\  a(\alpha) >0}}k(\alpha)a(\alpha)[\cot(\langle \alpha,\phi(t)\rangle) - \cot(\langle \sigma_0(\alpha),\phi(t)\rangle)]
\end{align*}
where $R_+^1 = \{e_i - e_j,\, 1 \leq i < j \leq m\}$ since $a(\alpha) = 0$ for $\alpha  \in \{e_i, 2e_i,\,i \neq m\}$, and $k(\alpha) = \beta$. 
Arguing as above, $T_{e_m} < \infty$ a.s. if $0< k_0 + k_1+1 < 2 \Leftrightarrow 0< p - (m-1) < 2/\beta$. \\
Finally, recall that $\mu := \pi/2 - \phi$ satisfies (\ref{JE1})  with $(p,q)$ intertwined, that is $2k_0 = \beta(q-p), k_1 = \beta(p-(m-1)) -1, 2k_2 = \beta$. Since $\mu_m = \pi/2 - \phi_1$ and since, as shown above, the drift term $F_t$ of $d\mu_m(t) = d\langle \mu(t),e_m\rangle$ is independent from $p,q$, we conclude that $\mu_m$ hits $0$ a.s. if 
$0 < \beta(q - (m-1)) < 2$ thereby $T_{\tilde{\alpha}} < \infty$ if $0 <  q - (m-1) < 2/\beta$.$\hfill \blacksquare$
\begin{nota}
When translated from $\phi$ to $\lambda$, the previous statements specialize to results from \cite{Dou}: in particular the eigenvalues process of both real and complex Jacobi processes never collide. However, it is known from \cite{Dou} that the eigenvalues process never hit $\{0,1\}$ when $p \wedge q \geq m -1 + 2/\beta$ for $\beta = 1,2$ respectively. This remains true for all $\beta > 0$ by the same arguments used in \cite{Dou}.   
\end{nota}

\section{Semi-group density}
We end this paper by writing down the semi group density of the $\beta$-Jacobi process $\lambda$. Before proceeding, we briefly consider two cases for which the semi group is known: the univariate case $m=1$ and the multivariate complex Hermitian one corresponding to $2$-Jacobi processes. In the former, the semi-group density reads (see \cite{Wong}) 
\begin{equation}\label{SGD1}
p_t^{r,s}(\theta,\lambda) = \sum_{n=0}^{\infty}e^{-2r_nt}\jac(\theta)\jac(\lambda) W^{r,s}(\lambda), \quad (\theta, \lambda) \in [0,1]^2,\, r,s > -1,
\end{equation} 
where $p = 2(r+1),q= 2(s+1), r_n = -n(n+r+s+1)$, $(\jac)_n$ are orthonormal Jacobi polynomials on $[0,1]$ and $W^{r,s}(\lambda) := C_{r,s} \lambda^r(1-\lambda)^{s}$ for some constant $C_{r,s}$ so that $W$ defines a probability distribution. No closed form is known for this density, nonetheless an attempt to get a handier expression was tried in \cite{Dem1}. 
In the latter, the $2$-Jacobi process is the Vandermonde transform of $m$ independent real Jacobi processes of parameters $d = 2(p-m+1): = 2(r+1), d'= 2(q-m+1):= 2(s+1)$ 
and conditioned never to collide (\cite{Dou}, p.141). More precisely, the Vandermonde function $V$ satisfies 
\begin{equation*}
\mathbf{L}V = -m(m-1)\left(\frac{2(m-2)}{3} + \frac{d+d'}{2}\right) V = -m(m-1)\left(\frac{2(m-2)}{3} + p+q - 2(m-1)\right) V:= cV 
\end{equation*}
where $\mathbf{L}$ is the generator of $m$ independent real Jacobi processes of parameters $d,d'$ (see appendix in \cite{Dou} p. 219). It follows by Karlin-McGregor formula (\cite{Kar}) that the semi group density is given by  
\begin{align*}
&p_t^{r,s,2}(\theta,\lambda) := e^{-ct}\frac{V(\lambda)}{V(\theta)} \det\left(\sum_{n=0}^{\infty}e^{-2n(n+r+s+1)t}\jac(\theta_i)\jac(\lambda_j)W^{r,s}(\lambda_j)\right)_{i,j=1}^m
\\& = e^{-ct}\det\left(\sum_{n=0}^{\infty}e^{-2n(n+r+s+1)t}\jac(\theta_i)\jac(\lambda_j)\right)_{i,j=1}^m\frac{W_m^{r,s,2}(\lambda)}{V(\theta)V(\lambda)}
\\&= e^{-ct}\left[\sum_{\sigma_1 \in S_m}\epsilon(\sigma_1)\sum_{n_1,\dots,n_m \geq 0}e^{-2\sum_{i=1}^mn_i(n_i + r + s+1)t}\prod_{i=1}^mP_{n_i}^{r,s}(\theta_i)P_{n_{i}}^{r,s} (\lambda_{\sigma_1(i)})\right]\frac{W_m^{r,s,2}(\lambda)}{V(\theta)V(\lambda)}
\\& = e^{-ct}\left[\sum_{\sigma_1,\sigma_2 \in S_m}\epsilon(\sigma_1)\sum_{n_1 \geq \dots \geq n_m \geq 0}e^{-2\sum_{i=1}^m n_{\sigma_2(i)}(n_{\sigma_2(i)} + r + s+1)t}
\prod_{i=1}^mP_{n_{\sigma_2(i)}}^{r,s}(\theta_i)P_{n_{\sigma_2(i)}}^{r,s}(\lambda_{\sigma_1(i)})\right]\frac{W_m^{r,s,2}(\lambda)}{V(\theta)V(\lambda)}
\end{align*}
where 
\begin{equation*}
W_m^{r,s,\beta}(\lambda) := C_{r,s,m,\beta} \prod_{i=1}^m\lambda_i^{r}(1-\lambda_i)^{s} \prod_{1\leq i < j\leq m}|\lambda_i - \lambda_j|^{\beta}
\end{equation*}   
$S_m$ is the symmetric group, and $\{0 < \lambda_m < \dots < \lambda_1 < 1\}, \{0 < \theta_m < \dots < \theta_1 < 1\}$. Note that, for a given partition $(n_1\geq \dots \geq n_m \geq 0)$ and $\sigma_2 \in S_m$, one has 
\begin{equation*}
\sum_{i=1}^m n_{\sigma_2(i)}(n_{\sigma_2(i)} + r + s+1) = \sum_{i=1}^mn_i(n_i + r + s+1).
\end{equation*}
Thus summing first over $\sigma_1$ with the change of variables $\sigma = \sigma_1\sigma_2$, one gets: 
\begin{align*}
p_t^{r,s,2}(\theta,\lambda) &= e^{-ct}\sum_{n_1 \geq \dots \geq n_m \geq 0}e^{-2\sum_{i=1}^mn_{i}(n_{i} + r + s+1)t} 
\frac{\det[P_{n_i}^{r,s}(\theta_j)]_{i,j=1}^m}{V(\theta)}\frac{\det[P_{n_i}^{r,s}(\lambda_j)]_{i,j=1}^m}{V(\lambda)}W_m^{r,s,2}(\lambda)
\\& = e^{-ct}\sum_{n_1 > \dots > n_m \geq 0}e^{-2\sum_{i=1}^mn_{i}(n_{i} + r + s+1)t} 
\frac{\det[P_{n_i}^{r,s}(\theta_j)]_{i,j=1}^m}{V(\theta)}\frac{\det[P_{n_i}^{r,s}(\lambda_j)]_{i,j=1}^m}{V(\lambda)}W_m^{r,s,2}(\lambda).
\end{align*} 
$p_t^{r,s,2}$ may be expressed as in (\ref{SGD1}) using the multivariate Jacobi polynomials $\jam$ defined by Lassalle which take a determinantal form for $\beta=2$ (\cite{Lass1}): 
\begin{equation*}
P_{\tau}^{r,s,2}(\lambda) = \frac{\det[P_{\tau_i+m-i}^{r,s}(\lambda_j)]_{i,j}}{V(\lambda)},
\end{equation*}
where $\tau = (\tau_1\geq \tau_2 \geq \cdots \tau_m)$. Thus, performing the index change $n_i := \tau_i + m-i$, one finally gets 
\begin{align}\label{SGD2}
p_t^{r,s,2}(\theta,\lambda) = \sum_{\tau_1 \geq \dots \tau_m \geq 0}e^{-2r_{\tau}^2} P_{\tau}^{r,s,2}(\theta)P_{\tau}^{r,s,2}(\lambda)W_m^{r,s,2}(\lambda),
\end{align}
where we set 
\begin{equation*}
r_{\tau}^2 = \sum_{i=1}^m\tau_i(\tau_i + r+s+1 + 2(m-i))
\end{equation*}
so that 
\begin{align*}
\sum_{i=1}^mn_{i}(n_{i} + r + s+1) + c/2 &= r_{\tau}^2 + \sum_{i=1}^m(m-i)(m-i+r+s+1) + c/2 
\\& =  r_{\tau}^2 + \sum_{i=1}^{m-1} i^2 + (p+q-2m+1)\sum_{i=1}^{m-1}i + c/2 
\\& = r_{\tau}^2  -\frac{(m-1)m(2m-1)}{3} + (p+q)\frac{m(m-1)}{2} + c/2
\\& = r_{\tau}^2.
\end{align*}
With regard to (\ref{SGD1}) and (\ref{SGD2}), it is natural to claim that : 
\begin{pro}
The semi group density of the $\beta$-Jacobi process is given on $\{0 < \lambda_m < \cdots < \lambda_1 <1\}$ by
\begin{equation}
\label{densite}
p_t^{r,s,\beta}(\theta, \lambda) := \sum_{\tau_1 \geq \cdots \geq \tau_m \geq 0}e^{-2r_{\tau}^{\beta}t}\jam(\theta)\jam(\lambda) W_m^{r,s}(\lambda)
\end{equation}
for $p \wedge q > (m-1) +1/\beta, $, where $\beta (p - (m-1)) := 2(r+1),\,\beta(q-(m-1)) := 2(s+1)$ and 
\begin{equation*}
r_{\tau}^{\beta} := \sum_{i=1}^m\tau_i(\tau_i +r+s+1 + \beta(m-i)).
\end{equation*}
\end{pro}

{\it Proof}: given a bounded symmetric function $f$ on $[0,1]^m$, define for $t > 0$
\begin{align*}
T_tf(\theta) &:= \int_{0 < \lambda_m< \dots < \lambda_1 < 1}f(\lambda)p_t^{r,s,\beta}(\theta, \lambda)d\lambda 
\end{align*}
for $\theta \in K := \{0 < \theta_1<\dots < \theta_m < 1\}$ and $T_0f = f$. The above integral converges: this uses the boundness of $f$, the exponential term with strictly positive $t$ and Fubini Theorem. Besides, $T_t{\bf 1} = 1$ and $||T_t||$ is bounded for all $t \geq 0$. Indeed, for $t > 0$, the first claim follows from Fubini's Theorem, the orthogonality of $(\jam)_{\tau}$ and $P_0^{r,s,\beta} = {\bf 1}$ so that the only non zero term is that corresponding to $\tau_i=0$ for all $i$. The proof of the second claim is almost similar. One also easily checks that $T_tT_s = T_{t+s}$. Let $\mathscr{L}_{\theta}$ be the generator of the $\beta$-Jacobi process $\theta$. Then, $\jam$ is an eigenfunction of $\mathscr{L}_{\theta}$  with eigenvalue $2r_{\tau}^{\beta}$ (\cite{Lass1}) so that $\mathscr{L}_{\theta}T_tf(\theta) = \partial_tT_tf(\theta)$ by the dominated convergence Theorem. 
Now, let us consider the Cauchy problem associated to $\mathscr{L}_{\theta}$: 
\begin{equation*}
\left\{\begin{array}{l}
 \partial_t u_{f}(t,\theta) = \mathscr{L}_{\theta}u_{f}(t,\theta) \\ 
 u_{f}(0,\cdot) = f, 
\end{array}\right.
\end{equation*}
where $u_{f} \in C^{1,2}(\re_+^{\star} \times K)  \cap C_b(\re^+ \cap K)$ with  boundary condition: 
\begin{equation*}
< \nabla u_f(t,\theta), n(\theta)> = 0  
\end{equation*}
where $n(\theta)$ is a unitary inward normal vector at $\theta \in \partial K$. Define $u_t (f)(\theta) := u_{f}(t,\theta)$. It is known (\cite{Str}) that the above Cauchy problem has a unique solution. Consequently, $u_t = T_t$ for all $t \geq 0$ and $(T_t)_{t \geq 0}$ is the semi group of the $\beta$-Jacobi process with density given by (\ref{densite}) (note that $T_tf \rightarrow  f$ in $L^2$-sense as $t \rightarrow 0$). $\hfill \blacksquare$ 

\section{Appendix}
Below, we give a quick proof of the fact that $\Delta h  > ch$ for some strictly negative constant $c$. To proceed, recall that, for $\phi \in A_{\circ}$, 
\begin{align*}
h(\phi) &= \prod_{i=1}^m\sin(\phi_i)\sin(2\phi_i)\prod_{i < j}\sin(\phi_i - \phi_j)\sin(\phi_i+\phi_j) 
\\& = 2^m\prod_{i=1}\sin^2(\phi_i)\cos(\phi_i)\prod_{j< i}(\sin^2(\phi_i) - \sin^2(\phi_j)) = 2^m \prod_{i=1}^m\lambda_i\sqrt{1-\lambda_i}V(\lambda) 
\end{align*}
where $\lambda_i = \sin^2(\phi_i)$ for $1 \leq i \leq m, \lambda = (\lambda_i)_{i=1}^m$ and $V$ stands for the Vandermonde function. Set 
\begin{equation*}
u(\lambda) = \prod_{i=1}^m\lambda_i\sqrt{1-\lambda_i},\quad v(\phi) = (\sin^2(\phi_1),\cdots, \sin^2(\phi_m)) = \lambda
\end{equation*}
so that $h(\phi) = 2^m(uV)(\lambda) = 2^m(uV)(v(\phi))$. It follows that 
\begin{align*}
\partial_i h(\phi) &= 2^m [2\sin(\phi_i)\cos(\phi_i)] \partial_i(uV)(\lambda)  = 2^m[2\sqrt{\lambda_i(1-\lambda_i)}] \partial_i(uV)(\lambda)
\end{align*}
and that 
\begin{align*}
\partial_i^2 h(\phi) &=  2^m[2\sqrt{\lambda_i(1-\lambda_i)}]\partial_i [2\sqrt{\lambda_i(1-\lambda_i)}\partial_i(uV)(\lambda)].
\\& = 2^m [4\lambda_i(1-\lambda_i)\partial_i^2(uV)(\lambda) + 2(1-2\lambda_i)\partial_i(uV)(\lambda)].
\end{align*}
Now, write 
\begin{equation*}
\partial_iu(\lambda) = \prod_{j\neq i}\lambda_j\sqrt{1-\lambda_j}\frac{2-3\lambda_i}{2\sqrt{1-\lambda_i}} = \frac{2-3\lambda_i}{2\lambda_i(1-\lambda_i)}u(\lambda)
\end{equation*}
and $\partial_i V = V\partial_i \log V$ so that 
\begin{equation*}
\partial_i (uV)(\lambda) = \left[\frac{2-3\lambda_i}{2\lambda_i(1-\lambda_i)} + \partial_i\log V(\lambda)\right] (uV)(\lambda) =  
\left[\frac{1}{\lambda_i} - \frac{1}{2(1-\lambda_i)} + \partial_i\log V(\lambda)\right] (uV)(\lambda)
\end{equation*}
and  
\begin{align*}
\frac{\partial_i^2 (uV)}{(uV)}(\lambda) &
= \left[\frac{1}{\lambda_i} - \frac{1}{2(1-\lambda_i)}+ \partial_i\log V(\lambda)\right]^2 + \left[\partial_i^2\log V(\lambda) - \frac{1}{\lambda_i^2}- \frac{1}{2(1-\lambda_i)^2}\right]. 
\\& = - \frac{1}{4(1-\lambda_i)^2} - \frac{1}{\lambda_i(1-\lambda_i)}+ [(\partial_i \log V(\lambda))^2 + \partial_i^2\log V(\lambda)] + 
\frac{2-3\lambda_i}{\lambda_i(1-\lambda_i)} \partial_i \log V(\lambda).
\end{align*}
Thus
\begin{align*}
\frac{\partial_i^2h}{h}(\phi) &= \frac{(1-2\lambda_i)(2-3\lambda_i)}{\lambda_i(1-\lambda_i)} + [2(1-2\lambda_i) + 4(2-3\lambda_i)]\partial_i \log V(\lambda) - 4 - 
\frac{\lambda_i}{1-\lambda_i} + 4\lambda_i(1-\lambda_i)\frac{\partial_i^2 V}{V}(\lambda)
\\& = (10-16\lambda_i)\partial_i \log V(\lambda) - 4 + \frac{2-7\lambda_i +5\lambda_i^2}{\lambda_i(1-\lambda_i)} + 4\lambda_i(1-\lambda_i)\frac{\partial_i^2 V}{V}(\lambda).
\end{align*}
Noting that $2-7\lambda_i +5\lambda_i^2 = (\lambda_i-1)(5\lambda_i -2)$ and using
\begin{equation*}
\sum_i\partial_i \log V(\lambda) = \sum_{j \neq i}\frac{1}{\lambda_i - \lambda_j} = 0, \quad \sum_{i\neq j} \frac{\lambda_i}{\lambda_i-\lambda_j} = \frac{m(m-1)}{2}
\end{equation*}
then
\begin{equation*}
\frac{\Delta h}{h}(\phi) = -8m(m-1) -4m + \frac{4}{V(\lambda)} \sum_{i=1}^m \lambda_i(1-\lambda_i)\partial_i^2 V(\lambda) - \sum_{i=1}^m\frac{5\lambda_i-2}{\lambda_i}.
\end{equation*}
Finally, from p. 219 in \cite{Dou}, 
\begin{equation*}
2\sum_{i=1}^m \lambda_i(1-\lambda_i)\partial_i^2 V(\lambda) = -\frac{2m(m-1)(m-2)}{3}V(\lambda)
\end{equation*}
so that
\begin{equation*}
\frac{\Delta h}{h}(\phi) = -8m(m-1) -9m  -\frac{4m(m-1)(m-2)}{3} +  \sum_{i=1}^m\frac{2}{\lambda_i} > -8m^2  - m  -\frac{4m(m-1)(m-2)}{3} := c
\end{equation*}
which finishes the proof.

\begin{nota}
In \cite{Gor}, the author comes to the $2$-Jacobi process as a limit in distribution of a rescaled Markov chain on the so-called Gelfand-Tsetlin graph.
\end{nota}


\begin{thebibliography}{99}
\bibitem{Bee}\emph{R. J. Beerends, E. M. Opdam}. Certain hypergeometric series related to the root system $BC$. {\it Trans. Amer. Math. Soc}.  {\bf 339},  no. 2. 1993, 581-607.
\bibitem{Bia}\emph{P. Biane}. Free Brownian motion, free stochastic calculus and random matrices. {\it Fie. Inst. Comm.}  {\bf 12}, Amer. Math. Soc. Providence, RI,  1997.
\bibitem{Cepa}\emph{E. C\'epa}. Equations  diff\'erentielles stochastiques multivoques. {\it S\'em. Proba.} {\bf XXIX}. 1995, 86-107. 
\bibitem{Cepa1}\emph{E. C\'epa, D. L\'epingle}. Diffusing particles with electrostatic repulsions. {\it Probab. Theo. Relat. Fields}. {\bf 107}, 1997, 429-449.
\bibitem{Cepa2}\emph{E. C\'epa, D. L\'epingle}. Brownian particles with electrostatic repulsion on the circle : Dyson's model for unitary random matrices revisited. {\it E. S. A. I. M : Probability and Statistics}. {\bf 5}. 2001, 203-224. 
\bibitem{Che}\emph{I. Cherednik, P. J. Forrester, D. Uglov}. Random matrices, Log-gases and the Calogero-Sutherland model. { \it Quantum Many-body Problems and Representation Theory}. MSJ Memoirs. 1998, 97-181.
\bibitem{Chy}\emph{O. Chybiryakov, N. Demni, L. Gallardo, M. Voit, M. R\"osler, M. Yor}. Harmonic and Stochastic Analysis of Dunkl operators. {\it Hermann Paris}, 2008.
\bibitem{De}\emph{N. Demni}. Generalized Bessel function of type $D$. {\it SIGMA. Special Volume: Dunkl operators and related topics}. {\bf 4}. 2008. 7 pages. 
\bibitem{Dem}\emph{N. Demni}. Radial Dunkl processes: existence, uniqueness and hitting time. {\it Submitted to C. R. A. S}.
\bibitem{Dem1}\emph{N. Demni, M. Zani}. Large deviations for statistics of Jacobi processes. {\it Stoch. Proc. Appl.} {\bf 119}, 2009. 518-533.
\bibitem{Dou}\emph{Y. Doumerc}. matrices al\'eatoires, processus stochastiques et groupes de r\'efl\'exions. {\it Ph. D. Thesis}, Paul Sabatier Univ. May 2005.
\bibitem{Gor}\emph{V. Gorin}. Non-colliding Jacobi processes as limits of Markov chains on Gelfand-Tsetlin graph. {\it To appear in J. Math. Sciences}. 
\bibitem{Gra}\emph{D. J. Grabiner}. Brownian motion in a Weyl chamber, non-colliding particles and random matrices. {\it Ann. I. H. P}. {\bf 35}, 1999, no. 2, 177-204. 
\bibitem{Hob}\emph{D. Hobson, W. Werner}. Non-colliding Brownian motions on the circle. {\it Bull. London Math. Society}. {\bf 28}. 1996, 643-650.
\bibitem{Hum}\emph{J. E. Humphreys}. Reflections Groups and Coxeter Groups. {\it Cambridge University Press}. {\bf 29}. 2000.
\bibitem{Kar}\emph{I. Karatzas, S. E. Shreve}. Brownian Motion and Stochastic Calculus, 2nd edition. {\it Springer-Verlag, New York}. 1991.
\bibitem{Karl}\emph{S. P. Karlin, G. McGregor}. Coincidence probabilities. {\it Pacif. J. Math}. {\bf 9}, 1959, 1141-1164.
\bibitem{Kat}\emph{M. Katori, H. Tanemura}. Symmetry of matrix-valued processes and noncolliding diffusion particle systems. {\it J. Math. Phy}. {\bf 45}, 2004, 3058-3085.
\bibitem{Kil}\emph{R. Killip, I. Nenciu.} Matrix models for circular ensembles. {\it Int. Math. Res. Not.} {\bf 50}, 2004, 2665--2701. 
\bibitem{Lass1}\emph{M. Lassalle}. Polyn\^omes de Jacobi g\'en\'eralis\'es. {\it C. R. A. S. Paris}. {\bf 312}. S\'erie {\bf I}, 1991, 425-428.
\bibitem{Rev}\emph{D. Revuz, M. Yor}. Continuous Martingales And Brownian Motion, $3^{\textrm{rd}}$ ed, Springer, 1999.
\bibitem{Scha}\emph{B. Schapira.} The Heckman-Opdam Markov processes.  {\it Probab. Theo. Relat. Fields.} 38, (2007), 495--519.
\bibitem{Str}\emph{D. W. Stroock, S.S. Varadhan.} Multidimensional Diffusion Processes. {\it Springer-Verlag}, 1971.
\bibitem{War}\emph{J. Warren, M. Yor}. The Brownian Burglar : Conditioning Brownian motion by its local time process. {\it S\'em. Proba.} {\bf XXXII.}, 1998, 328-342.
\bibitem{Wong}\emph{E. Wong}. The construction of a class of stationary Markov. {\it Proceedings. The $16^{th}$ Symposium. Applied Math}. AMS. Providence. RI. 1964. 264-276.
\end{thebibliography}
\end{document}